\documentclass[letterpaper, 10 pt, conference]{ieeeconf}  

\IEEEoverridecommandlockouts                              
\usepackage{graphics} 
\usepackage{booktabs}
\usepackage{epsfig} 
\usepackage{amsmath} 
\usepackage{amssymb}  
\usepackage{mathtools}
\let\temp\rmdefault
\usepackage{mathpazo}
\let\rmdefault\temp
\usepackage{tikz}
\usepackage{standalone}
\usepackage{bondgraphs}
\usepackage{url}
\usepackage{blox}
\usepackage{gensymb}
\usepackage{tabularx}
\usepackage{amsfonts,dcolumn}
\usepackage{blindtext}
\usepackage{rotating}
\usepackage[miama]{emf}
\DeclareMathAlphabet{\mathcal}{OMS}{cmsy}{m}{n}
\usetikzlibrary{shapes,arrows}
\usetikzlibrary{calc, patterns, decorations.pathmorphing, decorations.markings, shapes.geometric, arrows, positioning, angles}
\usepackage{hyperref}

\usepackage{textcomp}
\allowdisplaybreaks
\usepackage{accents}
\newcommand{\ubar}[1]{\underaccent{\bar}{#1}}
\newcommand{\R}{\mathbb{R}}
\DeclarePairedDelimiter{\norm}{\lVert}{\rVert}
\DeclarePairedDelimiter{\abs}{\lvert}{\rvert}%

\newtheorem{problem}{Problem}
\newtheorem{theorem}{Theorem}

\newtheorem{lemma}{Lemma}
\newtheorem{definition}{Definition}
\newtheorem{proposition}{Proposition}

\title{\LARGE \bf
Robust Estimator-Based Safety Verification: A Vector Norm Approach
}

\author{Binghan He$^{1}$, Gray C. Thomas and Luis Sentis
\thanks{
This work was supported by the U.S. Government and NASA Space Technology Research Fellowship NNX15AQ33H. 
We thank the members of the Human Centered Robotics Lab, University of Texas at Austin who provided insight and expertise that assisted the research. 
Authors are with The Departments of Mechanical Engineering (B.H., G.C.T.) and Aerospace Engineering (L.S.), University of Texas at Austin, Austin, TX. 
Send correspondence to $^{1}$$\;${\tt\small binghan at utexas dot edu}.
}
}

\newcommand\copyrighttext{%
\scriptsize 
Accepted for publication in American Control Conference (ACC)
\textcopyright 2020 IEEE. Personal use of this material is permitted. Permission from IEEE must be obtained for all other uses, in any current or future media, including reprinting/republishing this material for advertising or promotional purposes, creating new collective works, for resale or redistribution to servers or lists, or reuse of any copyrighted component of this work in other works.
DOI: \href{https://ieeexplore.ieee.org/document/9147276}{10.23919/ACC45564.2020.9147276}
}
\newcommand\copyrightnotice{%
\begin{tikzpicture}[remember picture,overlay]
\node[anchor=south,yshift=10pt] at (current page.south)
{\fbox{\parbox{\dimexpr\textwidth-\fboxsep-\fboxrule\relax}{\copyrighttext}}};
\end{tikzpicture}%
}

\usepackage{dcolumn}
\usepackage{multirow}
\begin{document}
\newcolumntype{L}[1]{>{\raggedright\arraybackslash}p{#1}}
\newcolumntype{C}[1]{>{\centering\arraybackslash}p{#1}}
\newcolumntype{R}[1]{>{\raggedleft\arraybackslash}p{#1}}

\maketitle
\thispagestyle{empty}
\pagestyle{empty}
\copyrightnotice

\vspace{-10pt}
\begin{abstract}
In this paper, we consider the problem of verifying safety constraint satisfaction for single-input single-output systems with uncertain transfer function coefficients. 
We propose a new type of barrier function based on a vector norm. This type of barrier function has a measurable upper bound without full state availability. 
An identifier-based estimator allows an exact bound for the uncertainty-based component of the barrier function estimate.
Assuming that the system is safe initially allows an exponentially decreasing bound on the error due to the estimator transient. 
Barrier function and estimator synthesis is proposed as two  convex sub-problems, exploiting linear matrix inequalities.
The barrier function controller combination is then used to construct a safety backup controller.
And we demonstrate the system in a simulation of a 1 degree-of-freedom human--exoskeleton interaction.
\end{abstract}

\section{Introduction}

Safe control is mission critical for robotic systems with humans in the loop. Uncertain robot model parameters and the lack of direct human state knowledge bring extra difficulty to the stabilization of human--robot systems. Methods such as robust loop shaping \cite{BuergerHogan2007TRO, HeThomasPaineSentis2019ACC, ThomasCoholichSentis2019AIM}, model reference adaptive control \cite{ChenChenYaoZhuZhuWangSong2016TMech} and energy shaping control \cite{lvGregg2017TCST} aim to balance the closed loop stability and performance of physical human robot interaction systems. However, there is no backup controller if these systems fail to maintain safety, because backup safety controllers require full state availability.

For systems with direct state measurements, safety is usually verified by a barrier certificate. Similar to a Lyapunov function, a barrier function or barrier certificate decreases at the boundary of its zero level set \cite{PrajnaJadbabaie2004HSCC}. While barrier certificate can be synthesized automatically through sum-of-squares (SoS) optimization \cite{Prajna2006Automatica}, a more ambitious goal is to combine the synthesis of the barrier function and the controller together through a control barrier function \cite{WielandAllgower2007IFAC}. Various methods such as backstepping \cite{TeeGeTay2009Automatica} and quadratic programming \cite{AmesXuGrizzleTabuada2016TAC, NguyenSreenath2016IJRR} create control barrier functions to ensure output and state constraint satisfaction while other methods such as semidefinite programming \cite{PylorofBakolas2016ACC} aimed to also include input saturation.

Safety warranties can also be considered a problem of finding an invariant set of the system which is also a subset the safe region in the state space. This allows us to consider using the synthesis of a quadratic Lyapunov function subject to the state and input constraints in a series of linear matrix inequalities (LMIs) \cite{BoydElGhaouiFeron1994Book}. To certify a larger safe region, composite quadratic Lyapunov functions can combine multiple existing certificates, either centered at the origin \cite{HuLin2003TAC} or with multiple equilibrium points \cite{ThomasHeSentis2018ACC}. The LQR-Tree strategy \cite{TedrakeManchesterTobenkinRoberts2010IJRR}, which could potentially be applied to safety control, creates a series of connected regions of attraction (also known as funnels) using quadratic Lyapunov functions for mapping the reachable state space. In \cite{HannafordRyu2002TRA}, a strategy was proposed to observe the safety of a system through its passivity which can be considered as a more conservative safety constraint than quadratic Lyapunov stability. 

A state space realization models a physical process if it correctly reproduces the corresponding output for each admissible input \cite{Morse1974CDC}. A Luenberger observer \cite{Luenberger1964TME} asymptotically estimates the state of such a model of a linear system with only the direct measurement of input and output. This idea has also been extended for system with nonlinear modeling error \cite{Zeitz1987SCL}. For bounded modeling errors, the estimation error converges to a residue set instead of zero \cite{CorlessTu1998Automatica}. Recently, a method of using sum-of-squares programming \cite{PylorofBakolasChan2019CSSL} aims to optimize the converging rate of a robust state estimation for uncertain nonlinear systems. But the estimated state still cannot be directly used for evaluation of barrier functions until it fully converges. The system could possibly violate the safety constraints before the barrier function estimation becomes valid.

In this paper, we aim to close the gap between state estimation and safety assurance for uncertain systems.
In order to address the barrier function estimation, we start with an identifier-based state estimator \cite{Morse1996TAC} which provides us a state estimate that is linear with the uncertain transfer function coefficients. Then, we define a vector norm based on a quadratic Lyapunov function such that a triangle inequality can be applied to decompose it into estimated state and estimation error. A convex polytopic bound on the estimated state is availiable through the estimator structure, and an upper bound on the estimation error arises from the convergence rate of the estimator and initial error. To obtain a larger safe (state-space) region, we extend this upper bound searching strategy to another vector norm defined based on a composite quadratic Lyapunov function \cite{HuLin2003TAC}, whose unit level set is a convex hull of the unit level sets of multiple quadratic Lyapunov functions. Using these vector norms, we derive our proposed barrier functions for uncertain systems with stable static output feedback. The synthesis of an estimator for the proposed barrier functions can be done in a two-step convex optimization using linear matrix inequalities, first optimizing the barrier function and then optimizing the estimator. This establishes a barrier pair \cite{ThomasHeSentis2018ACC}, which can be used with a hybrid safety controller to guarantee safety even for arbitrary inputs. In the end, our hybrid safety controller is demonstrated in a simulation of a simple human-exoskeleton interaction model with human stiffness uncertainty and velocity and force limits.

\section{Preliminaries}

\subsection{Problem Statement}

Let us consider an $n$-th order strictly proper uncertain SISO system $\Sigma_p$ with transfer function
\begin{flalign}
P (s) & = \frac{y(s)}{u(s)} = \frac{\qquad \, b_{1} s ^ {n-1} + \cdots + b_{n-1} s + b_{n}}{s ^ n + a_{1} s ^ {n-1} + \cdots + a_{n-1} s + a_{n}}, \label{tf} \\
a_i & \in [\ubar{a}_i, \, \bar{a}_i], \; \, i \in \{ 1, \, 2, \, \cdots , \, n \}, \label{ai} \\
b_j & \in [\ubar{b}_j, \, \bar{b}_j], \; \, j \in \{ 1, \, 2, \, \cdots , \, n \}, \label{bj}
\end{flalign}
where $u$ and $y$ are the input and output of $\Sigma_p$ and $i$ and $j$ are the indices of the polynomial coefficients. 

A state space realization of \eqref{tf} can be expressed as 
\begin{flalign}
\dot{x} & = A x + b_u u, \label{xp} \\ 
y &= c_0 x, \label{y} 
\end{flalign}
where $x$ is the state vector. We specify $(A, \, b_u, \, c_0)$ as an n-dimensional observable canonical form with $c_0 \overset{\Delta}{=} [1, \, 0, \, \cdots \, , 0]$. With the state space realization in the form of \eqref{xp}, the problem we consider is defined as follows.

\begin{problem}
Suppose there is exists a stable controller for the parameter uncertain system $\Sigma_p$ which can satisfy constraints $x\in \mathcal X$ and $u\in\mathcal U$ indefinitely for all initial states in $\mathcal X_s\subseteq \mathcal X$, find an estimator $\Sigma_e$ that can observe whether the system is inside the safe region $\mathcal X_s$ with direct measurement of only the input $u$ and output $y$ even when this controller is not necessarily active.
\end{problem}

\subsection{State Estimation}

Since only $u$ and $y$ are directly measured, we need to estimate $x$ in \eqref{xp} to verify safety. According to Lemma 1 in \cite{Morse1980TAC}, we can select a strictly stable $A_0$ in observable canonical form such that \eqref{xp} becomes
\begin{equation} \label{ss}
\begin{aligned} 
\dot{x} & = A_0 x + b_y y + b_u u,
\end{aligned}
\end{equation}
where $A$ in \eqref{xp} is replaced by $A_0 + b_y c_0$. Let the characteristic equation of $A_0$ be $s ^ n + \hat{a}_{1} s ^ {n-1} + \cdots + \hat{a}_{n-1} s + \hat{a}_{n}$. Since $(c_0, \, A_0)$ is also a pair in the observable canonical form, $b_y$ and $b_u$ are
\begin{flalign}
b_y & = [\hat{a}_1 - a_1, \, \hat{a}_2 - a_2, \, \cdots, \, \hat{a}_{n} - a_{n}] ^ T, \\
b_u & = [b_1, \, b_2, \, \cdots, \, b_{n}] ^ T,
\end{flalign}
which are either linear with or affine to the coefficients of the polynomials of $P(s)$ in \eqref{tf}.


We estimate $x$ through an identifier-based estimator which includes a pair of sensitivity function filters expressed as
\begin{equation} \label{id}
\begin{aligned} 
\dot{\theta}_{y} & = A_0 ^ T \theta_{y} + c_0 ^ T y, \\
\dot{\theta}_{u} & = A_0 ^ T \theta_{u} + c_0 ^ T u,
\end{aligned}
\end{equation}
where $(A_0 ^ T, c_0 ^ T)$ is a controllable pair in the canonical form.

\begin{lemma}
Suppose $E_y = C_0 ^ {-1} \Theta_{y} ^ T$ and $E_u = C_0 ^ {-1} \Theta_{u} ^ T$ where $C_0$ is the observability matrix of $(c_0, \, A_0)$, and $\Theta_{y}$ and $\Theta_{u}$ are the controllability matrices of $(A_0 ^ T, \, \theta_y)$ and $(A_0 ^ T, \, \theta_u)$. $E_y b_y + E_u b_u$ converges to $x$ exponentially.
\end{lemma}

\begin{proof}
This is similar to Lemma 2 in \cite{Morse1980TAC}. Notice that $C_0$ is also the transpose of the controllability matrix of $(A_0 ^ T, \, c_0 ^ T)$. We can derive from \eqref{id} that 
\begin{equation} \label{L-1-1}
\begin{aligned} 
\dot{E}_{y} ^ T & = A_0 ^ T E_{y} ^ T + I y, \\
\dot{E}_{u} ^ T & = A_0 ^ T E_{u} ^ T + I u.
\end{aligned}
\end{equation}
Because $A_0$ is in a canonical form, it is easy to show that $E_y A_0 = A_0 E_y$ and $E_u A_0 = A_0 E_u$. Therefore, by taking the transpose of \eqref{L-1-1}, we obtain $\dot{E}_{y} = A_0 E_{y} + I y$ and $\dot{E}_{u} = A_0 E_{u} + I u$.

If we define $\hat{x} \overset{\Delta}{=} E_y b_y + E_u b_u$, then $\dot{\hat{x}} = A_0 \hat{x} + b_y y + b_u u$.
Since $A_0$ is strictly stable, we have $x = \hat{x} + \epsilon$ where $\epsilon = e^{A_0 t} (x(0) - \hat{x} (0))$. 
\end{proof}

Notice that the dynamics of $\hat{x}$ can also be expressed as
\begin{equation} \label{luenberger}
\begin{aligned} 
\dot{\hat{x}} & = A \hat{x} + b_y (y - c_0 \hat{x}) + b_u u,
\end{aligned}
\end{equation}
which is a Luenberger observer of \eqref{xp}. However \eqref{luenberger} cannot be directly implemented because of the uncertainty in $b_y$ and $b_u$. The identifier-based estimator in \eqref{id} provides us a convex hull containing the estimated state vector $\hat x$,
\begin{equation} \label{xe_convex}
\begin{aligned}
\hat{x} (b_y, \, b_u) \in \text{Co} 
& \bigg\{
E_y
\begin{bmatrix*} 
\hat{a}_1 - a_1 \\
\hat{a}_2 - a_2 \\
\vdots \\
\hat{a}_n - a_n
\end{bmatrix*}
+
E_u
\begin{bmatrix*} 
b_1 \\
b_2 \\
\vdots \\
b_n
\end{bmatrix*}, \\
&\quad a_i \in \{\ubar{a}_i, \; \, \bar{a}_i\}, \, b_j \in \{\ubar{b}_j, \; \, \bar{b}_j\}, \\
& \quad \text{for} \; i, \, j = 1, \, 2, \, \cdots, \, n  \ \ \bigg\}.
\end{aligned}
\end{equation}
Because of the initial estimation error $\epsilon_0 \overset{\Delta}{=} x(0) - \hat{x}(0)$, any barrier function $B(x)$ aiming to constrain the system inside the safe region $\mathcal X_s$ cannot be directly bounded using $\hat{x} (b_y, \, b_u)$. Instead, our goal is to find an upper bound for the barrier function using both $\hat{x} (b_y, \, b_u)$ and $\epsilon_0$.

\subsection{Vector Norm Function}
In order to upper bound the barrier function proposed later in this paper, we recall the following two properties of a vector norm function.
\begin{lemma} \label{norm}
For every vector $x$ in some vector space within $\R ^ n$, let $\norm{\cdot}$ be a scalar function of $x$ with the following properties. 
\begin{itemize}
\item[(a)] $0 < \norm{x} < \infty$ except for $\norm{x} = 0$ at the origin.
\item[(b)] $\norm{\lambda x} = \abs{\lambda} \norm{x}$ for all $\lambda \in \R$.
\end{itemize}
Then $\norm{\cdot}$ satisfies 
\begin{itemize}
\item[(c)] $\norm{x + y} \leq \norm{x} + \norm{y}$
\end{itemize}
if and only if $\Omega \overset{\Delta}{=} \{x \mid \norm{x} \leq 1 \}$ is convex.
\end{lemma}

These properties (a), (b) and (c) in Lemma \ref{norm} are also called the three characteristic properties of a vector norm. 

\begin{lemma} \label{convex}
Let $\norm{\cdot}$ be a vector norm function satisfying (a), (b) and (c) in Lemma \ref{norm}. Suppose there is a vector $x_0 = \sum_{j = 1}^{N} \gamma_j x_j$ with $\sum_{j = 1}^{N} \gamma_j = 1$, $0 \leq \gamma_j <1$ for all $j = 1, \, 2, \, \cdots, \, N$ and $\norm{x_0} = \lambda$. Then there exists an index~$j$ such that $\norm{x_j} \geq \lambda$.
\end{lemma}

\begin{proof}
Suppose that $\norm{x_j} < \lambda$ for all $j = 1, \, 2, \, \cdots, \, N$. Based on (b) in Lemma \ref{norm}, we have $\norm{\gamma_j x} = \abs{\gamma_j} \norm{x}$ for all $j = 1, \, 2, \, \cdots, \, N$.

Applying (c) in Lemma~\ref{norm} to $\norm{x_0}$ we get 
\begin{equation} \label{conter}
\norm{x_0} = \big\|\sum_{j=1}^{N} \gamma_j x_j\big\| \leq \sum_{j=1}^{N} \gamma_j \norm{x_j} < \sum_{j=1}^{N} \gamma_j \lambda = \lambda,
\end{equation}
which contradicts $\norm{x_0} = \lambda$.
\end{proof}

\section{Barrier Estimation}

The triangle inequality of vector norms allows us to decompose the state $x$ into the estimated state $\hat{x}$ and estimation error $\epsilon$. While we do not know $x$, we know $\hat{x}$ and can bound $\epsilon$ within a decaying window---allowing us to extend barrier pairs \cite{ThomasHeSentis2018ACC} to systems without full state availability.

\subsection{Norm of Quadratic Lyapunov Function}

Let us define a quadratic Lyapunov function as $V_q (x) = x ^ T Q ^ {-1} x$ where $Q$ is a positive definite matrix. We can form a vector norm using its square root,
\begin{equation} \label{norm_q}
\norm{x}_q \overset{\Delta}{=} V_q ^ {\frac{1}{2}} (x),
\end{equation}
because $V_q (x)$ is positive definite, $V_q (\lambda x) = \lambda ^ 2 V_q (x)$ and the unit level set of $V (x)$ is convex.

For a given pair of $E_y$ and $E_u$, we can derive from Lemma~\ref{convex} that the maximum value of $\norm{\hat{x}}_q$ occurs at one of vertices of the convex hull in \eqref{xe_convex}.

\begin{theorem} \label{bound_vq}
If a strictly stable matrix $A_0$ in \eqref{ss} and \eqref{id} satisfies
\begin{equation} \label{lp}
A_0 Q + Q A_0 ^ T + 2 \alpha Q \preceq 0,
\end{equation}
then for all $t \geq 0$ there exists an $i \in \{ 1, \, \cdots, \, N \}$ such that
\begin{equation} \label{bp}
\norm{x}_q \leq \norm{\hat{x} (b_{yi}, \, b_{ui})}_q + e ^ {- \alpha t} \norm{\epsilon_0}_q,
\end{equation}
where $\hat{x} (b_{yi}, \, b_{ui})$ for $i = 1, \, 2, \, \cdots, \, N$ are the all vertices of \eqref{xe_convex}.
\end{theorem}

\begin{proof}
From Lemma \ref{norm}, we have $\norm{x}_q \leq \norm{\hat{x}}_q + \norm{\epsilon}_q$. 
The time derivative of $V_q (\epsilon)$ can be expressed as
\begin{equation} \label{p-3-2}
\dot{V}_q (\epsilon) \text{\small $= \epsilon ^ T (Q ^ {-1} A_0 + A_0 ^ T Q ^ {-1}) \epsilon = \epsilon ^ T Q ^ {-1} (A_0 Q + Q A_0 ^ T) Q ^ {-1} \epsilon$}.
\end{equation}
By substituting \eqref{lp}, $\dot{V}_q (\epsilon) \leq - 2 \alpha V_q (\epsilon)$
which guarantees that $V_q(\epsilon) \leq e ^ {- 2 \alpha t} V_q (\epsilon_0)$. Therefore, $\norm{\epsilon}_q \leq e ^ {- \alpha t} \norm{\epsilon_0}_q$. Together with Lemma \ref{convex}, we have \eqref{bp}.
\end{proof}

This Theorem~\ref{bound_vq} provides an upper bound on $\norm{x}_q$ which is available in that it be calculated from $\hat{x}$ and $\epsilon_0$ for all $t\geq0$. 

\subsection{Norm of Composite Quadratic Lyapunov Function}

\begin{figure}[!tbp]
    \footnotesize
    \centering
    	\def\svgwidth{.42\textwidth}
    	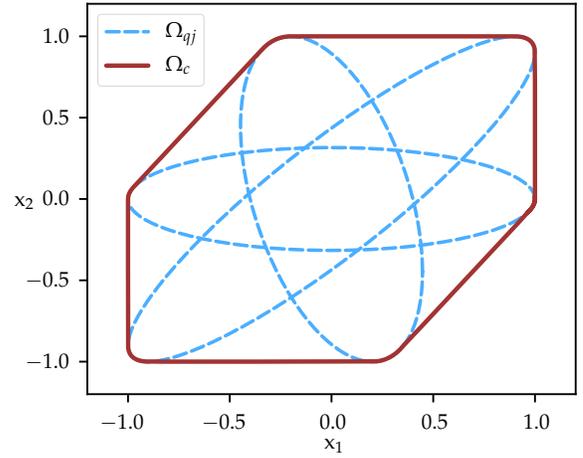
    \caption{A unit ball $\Omega_{c}$ of $\norm{x}_c$ equivalent to the convex hull of the ellipsoidal unit balls $\Omega_{qj}$ of three different $\norm{x}_{qj}$.}
    \vspace{-10pt}
    \label{fig:unit_ball}
\end{figure}

In order to obtain a larger safe region $\mathcal{X}_s$, a composite quadratic Lyapunov function is considered. For multiple different quadratic Lyapunov functions defined with positive-definite matrices $Q_1$, $Q_2$, $\cdots$, $Q_{n_q}$, a composite quadratic Lyapunov function \cite{HuLin2003TAC} is defined as 
\begin{align}
V_c (x)   & \overset{\Delta}{=} \min_{\gamma} \, x ^ T Q ^ {-1} (\gamma) x, \label{vc} \\
Q(\gamma) & \overset{\Delta}{=} \sum_{j = 1}^{n_q} \gamma_j Q_j, \label{qj}
\end{align}
where $\sum_{j = 1}^{n_q} \gamma_j = 1$ and $\gamma_j \geq 0$ for all $j = 1, \, 2, \, \cdots, \, n_q$. 

The unit level set of $V_c (x)$ is the convex hull of all the unit level sets of $V_{qj} (x) = x^T Q_j ^ {-1} x$ for $j = 1, \, 2, \, \cdots, \, n_q$ and is therefore also a convex shape (see Fig.~\ref{fig:unit_ball}). Since $V_c (x)$ is positive definite and $V_c (\lambda x) = \lambda ^ 2 V_c (x)$, we can use Lemma~\ref{norm} to define the ``composite'' vector norm,
\begin{equation} \label{norm_c}
\norm{x}_c \overset{\Delta}{=} V_c ^ {\frac{1}{2}} (x).
\end{equation}




\begin{theorem} \label{bound_vc}
For all $j = 1, \, 2, \, \cdots, \, n_q$, if a strictly stable matrix $A_0$ in \eqref{ss} and \eqref{id} satisfies
\begin{equation} \label{lc}
A_0 Q_j + Q_j A_0 ^ T + 2 \alpha Q_j \preceq 0, 
\end{equation}
then for all $ t\geq 0$ there exists an $i \in \{ 1, \, \cdots, \, N \}$ such that
\begin{equation} \label{bc}
\norm{x}_c \leq \norm{\hat{x} (b_{yi}, \, b_{ui})}_c + e ^ {- \alpha t} \norm{\epsilon_0}_c
\end{equation}
where $\hat{x} (b_{yi}, \, b_{ui})$ for $i = 1, \, 2, \, \cdots, \, N$ are the all vertices of \eqref{xe_convex}.
\end{theorem}

\begin{proof}
As in Theorem \ref{bound_vq}. 
\end{proof}

Therefore, an upper bound on $\norm{x}_c$ can be calculated using $\hat{x}$ and $\epsilon_0$ for all $t\geq0$.

\subsection{Barrier Pairs}

\begin{definition}[See \cite{ThomasHeSentis2018ACC}] \label{barrier_pair}
A \emph{Barrier Pair} is a pair of functions $(B,\ k)$ with two following properties:
\begin{itemize}
\item[(a)] $-1<B(x)\leq 0,\quad u=k(x) \implies \dot B(x) < 0$,
\item[(b)] $B(x)\leq 0 \implies x \in \mathcal X,\ k(x) \in \mathcal U$,
\end{itemize}
where (a) and (b) are also called invariance and constraint satisfaction.
\end{definition}

We can now introduce two barrier pairs using our vector norms and static output feedback controller.

\begin{proposition} \label{barrier_pair_q}
Suppose $V_q$ is a quadratic Lyapunov function for system of \eqref{xp} and \eqref{y} with a static output feedback $u = \mathbf k y$ and $\Omega_q$ is a unit ball of $\norm{x}_q$ defined as \eqref{norm_q}. If 
\begin{equation} \label{Omega_q}
\Omega_q \subseteq \mathcal{X} \cap \{x \mid c_0 x \in \mathbf k ^ {- 1} \mathcal{U} \},
\end{equation}
then $(\norm{x}_q-1, \ \mathbf k y)$ is a barrier pair.
\end{proposition}

\begin{proof}
Let $B_q (x) = \norm{x}_q - 1$. Its time derivative is
\begin{equation} \label{B_q_dot}
\dot{B}_q (x) = \frac{1}{2} \norm{x}_q ^ {-1} \dot{V}_q.
\end{equation}
Since $\norm{x}_q ^ {-1} > 0$ and $\dot{V}_q < 0$ when $-1< B_q (x) \leq 0$, $(B_q (x), \ \mathbf k y)$ satisfies (a) in Definition \ref{barrier_pair}. From \eqref{Omega_q}, we also have (b) in Definition \ref{barrier_pair} satisfied.
\end{proof}

\begin{proposition}
Suppose $V_c$ is a composite quadratic Lyapunov function defined as \eqref{vc} and \eqref{qj} for the system of \eqref{xp} and \eqref{y}, with static output feedback $u = \mathbf k y$ and that $\Omega_c$ is a unit ball of $\norm{x}_c$ defined as in \eqref{norm_c}. If we have
\begin{equation} \label{Omega_c}
\Omega_c \subseteq \mathcal{X} \cap \{x \mid c_0 x \in \mathbf k ^ {- 1} \mathcal{U} \},
\end{equation}
then $(\norm{x}_c-1,\ \mathbf k y)$ is a barrier pair.
\end{proposition}

\begin{proof} \label{barrier_pair_c}
As in Proposition \ref{barrier_pair_q}.
\end{proof}

As in \eqref{bp} and \eqref{bc}, upper bounds of the barrier functions of these two barrier pairs can be calculated using $\hat{x}$ and $\epsilon_0$ for all $t\geq0$.



\section{Synthesis}

Both barrier functions $B_c (x) \overset{\Delta}{=} \norm{x}_c - 1$ and their identifier-based estimators can be synthesized with LMIs, through the sub-problem of synthesizing $B_q (x) \overset{\Delta}{=} \norm{x}_q - 1$.




\subsection{\texorpdfstring{$Q_j$}{TEXT} Synthesis}

With static output feedback, the closed loop system of \eqref{xp} is still a polytopic linear differential inclusion (PLDI) model \cite{BoydElGhaouiFeron1994Book} $\dot{x} \in A_c x$ with
\begin{align} 
A_c = \text{Co} 
& \bigg\{
\begin{bmatrix*} 
0      & 1 &        & 0 \\
\vdots &   & \ddots &   \\
0      & 0 &        & 1 \\
0      & 0 & \cdots & 0
\end{bmatrix*}
-
\begin{bmatrix*} 
a_1 \\
a_2 \\
\vdots \\
a_n
\end{bmatrix*}
c_0 +
\begin{bmatrix*} 
b_1 \\
b_2 \\
\vdots \\
b_n
\end{bmatrix*}
\mathbf k c_0, \nonumber \\
& \quad a_i \in \{\ubar{a}_i, \; \, \bar{a}_i\}, \, b_j \in \{\ubar{b}_j, \; \, \bar{b}_j\}, \nonumber\\
& \quad\text{for} \; i, \, j = 1, \, 2, \, \cdots, \, n \ \  \bigg\}.
\label{LDI}
\end{align}
Supposing that $\mathcal{X}$ and $\mathcal{U}$ can be described (perhaps conservatively) as
\begin{align}
\mathcal X & = \{x : \, |f_i x| \leq 1, \ i = 1, \, 2, \, \cdots, \, n_f \}, \\
\mathcal U & = \{u :  \quad |u| \leq \bar{u} \}, 
\end{align}
they can be enforced by LMIs
\begin{align}
f_i Q_j f_i ^ T & \leq 1, \quad \forall\ i = 1, \, 2, \, \cdots, \, n_f, \label{eq:state_lim} \\
c_0 Q_j c_0 ^ T & \leq \frac{\bar{u} ^ 2}{\mathbf k ^ 2}. \label{eq:input_lim}
\end{align}
To synthesize $Q_j$, we maximize the width of the unit ball of $x^T Q_j ^ {-1} x$ along some state space direction $x_j$ by minimizing $\rho_j$ subject to the following LMI
\begin{equation}
\begin{bmatrix}
\rho_j & x_j^T\\
x_j & Q_j
\end{bmatrix}\succeq 0,\label{eq:drection_lim}
\end{equation}
such that the optimization sub-problem becomes
\begin{equation} \label{optimization_Qj}
\begin{aligned}
& \underset{Q_j}{\text{minimize}}
& & \rho_j \\
& \text{subject to} & & \eqref{eq:state_lim}, \ \eqref{eq:input_lim}, \ \eqref{eq:drection_lim},\ Q_j \succ 0, \\
&&& A_{ci} Q_j + Q_j A_{ci}^{T} + 2 \alpha_0 Q_j \preceq 0, \ \\
&&& \forall \ i = 1, \, 2, \, \cdots, \, N.
\end{aligned}
\end{equation}
where $A_{ci}$ for $i = 1, \, 2, \, \cdots, \, N$ are the all vertices of \eqref{LDI}. A positive value of $\alpha_0$ is used to guarantee a minimum exponential decay rate for $\norm{x}_{c}$.

\begin{figure}[!tbp]
  	\centering
    \includestandalone[width=.42\textwidth]{block-diagram}
    \caption{Block diagram consisting of plant $\Sigma_p$, estimator $\Sigma_e$ and hybrid safety controller $\Sigma_s$.}
    \vspace{-5pt}
    \label{fig:block-diagram}
\end{figure}

\begin{figure}
    \small
    \centering
    	\def\svgwidth{.28\textwidth}
    	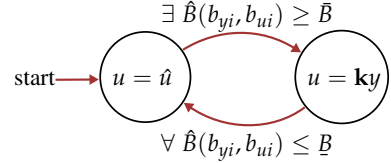
    \caption{Switching logic of hybrid safety controller $\Sigma_{s}$.}
    \vspace{-10pt}
    \label{fig:safety-controller}
\end{figure}

\subsection{\texorpdfstring{$A_0$}{TEXT}  Synthesis}
While it is simple to specify an $A_0$ in \eqref{ss} with a fast decay rate of $\epsilon$ (choosing big negative-real-part eigenvalues), this does not necessarily improve the value of $\alpha$ in \eqref{bc}. To synthesize an $A_0$ in the set of matrices in observable canonical form $\mathcal O\subset \R ^ {n \times n}$ we directly optimize for $\alpha$:
\begin{equation} \label{optimization_A0}
\begin{aligned}
& \underset{A_0\in\mathcal{O}}{\text{maximize}}
& & \alpha \\
& \text{subject to} & & A_0 Q_j + Q_j A_0^T + 2 \alpha Q_j \preceq 0, \\
&&& \forall \ j = 1, \, 2, \, \cdots, \, n_q,
\end{aligned}
\end{equation}
knowing that a solution $\alpha \geq \alpha_0$ will exist. (Any $A_0$ in the convex hull of \eqref{LDI} is guaranteed to satisfy the constraints in \eqref{optimization_A0} with a decay rate of $\alpha_0$.)

\subsection{Hybrid Safety Controller}

To enforce safety satisfaction on a potentially unsafe input $\hat{u}$, we can estimate the barrier function as 
\begin{equation} \label{be}
\hat{B}_c(b_{y}, \, b_{u}) \overset{\Delta}{=} \norm{\hat{x} (b_{y}, \, b_{u})}_c + e ^ {- \alpha t} \norm{\epsilon_0}_c - 1.
\end{equation}
With this estimate, we can design a hybrid safety controller $\Sigma_s$ which decides whether to apply either $\hat{u}$ or $\mathbf{k} y$ (that is, the safety backup control law) as the input in order to keep $B_c\leq 0$ (see Fig.~\ref{fig:block-diagram}) and therefore guarantee safety. 

According to Theorem \ref{bound_vc}, system $\Sigma_p$ is guaranteed to be safe if $\hat{B}_c(b_{yi}, \, b_{ui}) \leq 0$ for all vertices $\hat{x} (b_{yi}, \, b_{ui})$ of convex hull \eqref{xe_convex}. Therefore, the switching logic for $\Sigma_s$ defined in Fig.~\ref{fig:safety-controller}, which introduces two near-zero thresholds $\ubar{B}$ and $\bar{B}$ (with $-1 < \ubar{B} < \bar{B} \leq 0$), will result in robust safety.

\section{Example}

\begin{figure}[!tbp]
    \centering
    	\begin{tikzpicture}
\tikzstyle{spring}=[thick,decorate,decoration={zigzag,pre length=0.2cm,post length=0.2cm,segment length=6}]
\tikzstyle{damper}=[thick,decoration={markings,  
  mark connection node=dmp,
  mark=at position 0.5 with 
  {
    \node (dmp) [thick,inner sep=0pt,transform shape,rotate=-90,minimum width= 5pt,minimum height=1pt,draw=none] {};
    \draw [thick] ($(dmp.north east)+(2pt,0)$) -- (dmp.south east) -- (dmp.south west) -- ($(dmp.north west)+(2pt,0)$);
    \draw [thick] ($(dmp.north)+(0,-2pt)$) -- ($(dmp.north)+(0,2pt)$);
  }
}, decorate]
\tikzstyle{ground}=[fill,pattern=north east lines,draw=none,minimum width=0.75cm, minimum height=0.15cm]

\begin{scope}[xshift=5cm]
\node [draw,outer sep=0pt,thick] (M1) [minimum width=1.0cm, minimum height=0.6cm] {\footnotesize $m_e=1$};

\node (ground1) [ground,anchor=north, yshift=-0.25cm, minimum width=1.5cm] at (M1.south) {};

\draw (ground1.north east) -- (ground1.north west);

\draw (M1.south west) ++ (0.2cm,-0.125cm) circle (0.125cm)  (M1.south east) ++ (-0.2cm,-0.125cm) circle (0.125cm);

\node (wall1) [ground, rotate=-90, minimum width = 0.56cm, yshift = - 1.6cm] {};

\draw (wall1.north east) -- (wall1.north west);

\node (wall2) [ground, rotate=-90, minimum width = 0.56cm, yshift =   1.6cm] {};

\draw (wall2.south east) -- (wall2.south west);

\draw [damper] ($(M1.north east)!(wall2.90)!(M1.south east)$) -- ($(wall2.south east)!(wall2.90)!(wall2.south west)$) node [midway,below] {\footnotesize $b_e$};

\draw [spring] ($(wall1.north east)!(wall1.90)!(wall1.north west)$) -- ($(M1.north west)!(wall1.90)!(M1.south west)$) node [midway,below] {\footnotesize $k_h$};

\draw [-latex,ultra thick] ($(M1.north east)$) ++ (0.cm,-0.0cm) -- +( .5cm,0cm) node [midway,above] {\footnotesize $u$};

 \draw[thick, dashed] ($(M1.north west)$) -- ($(M1.north west) + (0,.5)$);
 \draw[thick, -latex] ($(M1.north west) + (0,0.25)$) -- ($(M1.north west) + (0.5,0.25)$) node [midway, above] {\footnotesize $\mathrm{x}_e$};
  \draw[thick, dashed] ($(wall1.north west)$) -- ($(wall1.north west) + (0,.5)$);
 \draw[thick, -latex] ($(wall1.north west) + (0,0.25)$) -- ($(wall1.north west) + (0.5,0.25)$) node [midway, above] {\footnotesize $\mathrm{x}_h$};

\end{scope}
\end{tikzpicture}
    \caption{Our simplified human--exoskeleton interaction model, a mass-spring-damper system, includes an uncertain human stiffness $k_h$, an exoskeleton damping $b_e$, and an exoskeleton inertia $m_e$.}
    \vspace{-5pt}
    \label{fig:mass-damper-spring}
\end{figure}
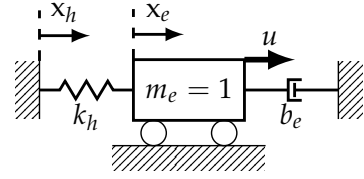

To illustrate robust barrier function estimation and hybrid safety control, we introduce a simplified human-exoskeleton interaction model. As shown in Fig.~\ref{fig:mass-damper-spring}, this model is a mass-spring-damper with uncertain human stiffness $k_h$, the exoskeleton damping $b_e$, and exoskeleton inertia $m_e$. 

\subsection{Human-Exoskeleton Interaction Model}
The exoskeleton plant can be expressed as a transfer function
\begin{equation}
P(s) = \frac{y (s)}{u (s)} = \frac{k_h}{m_e s ^ 2 + b_e s + k_h}
\end{equation}
where the input $u$ is the actuator force exerted and the output $y\overset{\Delta}{=} k_h (\mathrm{x}_e - \mathrm{x}_h)$ is the contact force between human and exoskeleton. 
Although the contact force and the exoskeleton position, $\mathrm{x}_e$, can be measured, the reference position, $\mathrm{x}_h$, of the human spring is not available because of the unknown stiffness.

Suppose that the uncertain value of $k_h$ is in the range from 4 to 12 and that the value of $b_e$ is 12. We can express the closed loop $A_c$ matrix set with static output feedback as a convex hull
\begin{equation} \label{LDI_example}
\begin{aligned}
A_c = \text{Co} \bigg\{ 
&
\begin{bmatrix*} 
- 12 & 1 \\
- k_h & 0 \\ 
\end{bmatrix*}
+
\begin{bmatrix*} 
0   \\
k_h \\
\end{bmatrix*}
\mathbf k c_0, 
& k_h = 4, \; 12 \bigg\}.
\end{aligned}
\end{equation}
And the safety constraints can be defined via the sets
\begin{equation}
\begin{aligned}
 \mathcal X &= \{ [ \mathrm{x}_1, \, \mathrm{x}_2 ] ^ T: \ |-\mathrm{x}_1 + \mathrm{x}_2/12| \leq 1 \},  \\
\mathcal U &= \{u :\ |u| \leq 1.2  \},
\end{aligned}
\end{equation}
where the output $y=\mathrm{x}_1$ and $\dot y=-12 \mathrm x_1+\mathrm x_2$, so $\mathcal X$ is constraining the output derivative $|\dot{y}| \leq 12$.

\subsection{Simulation}
\begin{figure*}[!tbp]
    \small
    \centering
    	\def\svgwidth{.99\textwidth}
    	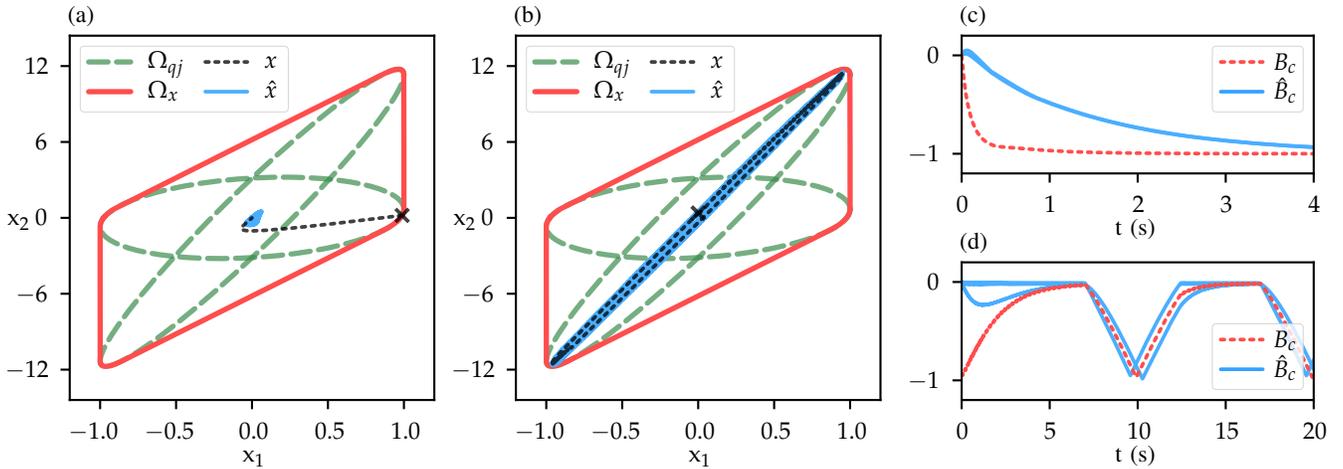
    \vspace{-9pt}
    \caption{Simulation results. In the fist simulation, the system state is initialized near the boundary of $\Omega_c$ (phase plot in (a)). The maximum $\hat{B}_c (k_h)$ converges slower than $B_c$ (in (c)). In the second simulation, an unsafe sinusoidal input $\hat{u}$ is forced to be safely inside $\Omega_c$ by a hybrid safety controller (phase plot in (b)). This safety controller activates only when $\max(\hat{B}_c (k_h))\approxeq 0$ (see (d)).%
    }
    \vspace{-12pt}
    \label{fig:simulation}
\end{figure*}

We choose the static output feedback gain $\mathbf{k}=-1.2$ which is stable, and leads to a human amplification factor of $2.2$ and the output constraint $|y| \leq 1$. In Fig.~\ref{fig:simulation}, we construct a barrier function $B_c$ from two different quadratic Lyapunov functions (optimized along directions $x_j=[1, \, 0]^T$ and $x_j=[1, \, 12]^T$) generated by synthesis \eqref{optimization_Qj} with a shared decay rate of $\alpha_0 = 0.50$. Then, an $A_0$ matrix with a characteristic polynomial $s ^ 2 + \hat{a}_1 s + \hat{a}_2$ (with $\hat{a}_1 = 13.60$ and $\hat{a}_2 = 18.68$) and a higher decay rate ($\alpha = 0.68$) is generated from synthesis \eqref{optimization_A0}. Notice that the optimal $A_0$ matrix is not exactly inside the convex hull of \eqref{LDI_example}.

In the numerical simulation, human stiffness $k_h=8$. In our first simulation, we release the system near the boundary of $\Omega_c$ with zero nominal input. In the second simulation the system starts at the origin and we apply a nominal input $\hat{u}$ which tracks an unsafe reference $y$ trajectory: $y(t) = 1.2\cdot \sin(0.05 \cdot 2 \pi t)$. In the first test (Fig.~\ref{fig:simulation}.a) the static output feedback is always on, to demonstrate the slower decay of $\hat{B}_c(k_h)$. In the second (Fig.~\ref{fig:simulation}.b), the static output feedback is turned on when $\max(\hat{B}_c (k_h))\geq \bar{B}=-0.01$ and is turned off when all values of $\max(\hat{B}_c (k_h))\leq\ubar{B}=-0.02$. This switching logic forces the system to stop near the boundary of $\Omega_c$---deviating from the unsafe trajectory to produce a safe output. In both tests, the largest element of $\hat{B}_c(k_h)$ converges to zero slower than the value of $B_c$ such that $\max(\hat{B}_c(k_h))\geq B_c$, as shown in Fig.~\ref{fig:simulation}.c and Fig.~\ref{fig:simulation}.d respectively. 

\setlength{\belowcaptionskip}{-12pt}

\section{Discussion}

Because the estimated barrier function $\hat{B}_c$ includes the transient term $e^{-\alpha t}\|\epsilon_0\|_c$, which is initially at a value of 1---indicating our assumption that the system state starts within the safe region---the safety backup controller is always active initially. Since our system is modelled as noiseless, this transient decays towards zero and the composite barrier approaches a steady state that robustifies the barrier only against the effect of parameter uncertainty. If there were additional uncertainty in the measurements, for example Gaussian noise, then this transient would asymptotically approach a steady-state non-zero value instead. 

When the barrier function estimate is near the thresholds, as in Fig.~\ref{fig:simulation}.d, the control switches rapidly between $\hat{u}$ and $\mathbf{k} y$. The frequency of this switching can be controlled by increasing the gap between $\ubar{B}$ and $\bar{B}$, which can reduce the risk of activating unmodeled high frequency dynamics. 


If we were to combine the synthesis of the barrier function and estimator together in one optimization, setting $\alpha$ equal to $\alpha_0$, the LMIs in \eqref{optimization_A0} would become bilinear matrix inequalities (BMIs). Since the estimator state equation in \eqref{luenberger} is in a canonical form, a new variable could be defined using the sufficient condition proposed in \cite{CrusiusTrofino1999TAC} to reduce these BMIs to LMIs. However, this sufficient condition would not always result in a feasible optimization problem.

The proposed synthesis strategy can be extended to systems with a pre-specified dynamic output feedback---albeit inefficiently. Though the states of a dynamic output feedback controller are perfectly known, they can also be considered part of the plant and estimated the same way as plant states. Input constraints can be incorporated as state constraints on the part of this composite plant that represents controller states. Then, synthesis \eqref{optimization_Qj} and \eqref{optimization_A0} can be applied directly by setting $\mathbf{k}$ to be zero.

This strategy can also be naturally extended to supervisory control \cite{Morse1996TAC} of a family of sub-optimal output feedback controllers subject to different subsets of parameter uncertainty by turning \eqref{id} into a shared state parameter identifier, as proposed in \cite{Morse1980TAC}. Substituting for methods such as constrained model reference adaptive control \cite{arabi2018set, lafflitto2018barrier} and adaptive control barrier functions \cite{taylor2019adaptive}, safety during the parameter adaptation could be enforced by switching to a backup barrier pair which is robust to the full parameter uncertainty.







\bibliographystyle{IEEEtran}
\bibliography{main}

\end{document}